\documentclass[a4paper,twoside,12pt]{article}
\usepackage{amssymb} 
\begin{document}
\title{A Conjecture on Zero-sum 3-magic Labeling of 5-regular Graphs \footnotemark[2]
\author{Guanghua Dong$^{1}$ and Ning Wang$^{2}$\\
{\small\em 1.Department of Mathematics, Tianjin Polytechnic
University, Tianjin, 300387, China}\\
\hspace{-5mm}{\small\em 2.Department of Information Science and
Technology,
Tianjin University of Finance }\\
\hspace{-74mm} {\small\em and Economics, Tianjin, 300222, China}\\
}} \footnotetext[2]{\footnotesize \em   This work was partially
Supported  by the China Postdoctoral Science Foundation funded
project (Grant No: 20110491248) and the National Natural Science
Foundation of China (Grant No: 11301381).}
\footnotetext[1]{\footnotesize \em E-mail: gh.dong@163.com(G. Dong);  ninglw@163.com(N. Wang). }

\date{}
\maketitle

\vspace{-.8cm}
\begin{abstract}

In this paper, we obtained that every 5-regular graph admits a zero-sum 3-magic labeling, which give an affirmative answer to a conjecture proposed by Saieed Akbari, Farhad Rahmati and Sanaz Zare in $Electron.$ $J.$ $Combin.$.


\end{abstract}


\bigskip
\noindent {\bf 1. Introduction}

Graph considered here are all finite and undirected with vertex set $V(G)$ and edge set $E(G)$. A $multigraph$ is a graph with multiple edges. If every vertex in a graph has the same degree $r$ then this graph is referred to as a $r$-$regular$ graph. A matching $M$ in $G$ is a set of independent edges, and $|M|$ denotes the number of edges in $M$. A $factor$  of a graph $G$ is a spanning subgraph of $G$. A $k$-$factor$ of $G$ is a factor of $G$ that is $k$-$regular$. Thus a $1$-$factor$ of $G$ is a matching that saturates all vertices of $G$, and  is called a $perfect$ $matching$ of $G$. A mapping $l : E(G) \rightarrow A$, where $A$ is an abelian group which written additively,  is called a $labeling$ of the graph $G$. Given a labeling $l$ of the graph $G$, the symbol $s(v)$, which represents the sum of the labels of  edges incident with $v$, is defined to be $s(v) = \sum_{uv\in E(G)} l(uv)$, where $v\in V(G)$. For every positive integer $h\geqslant 2$, a graph $G$ is said to be $zero$-$sum$ $h$-$magic$ if there is an edge labeling from $E(G)$ into $\mathbb{Z}_{h} \setminus \{0\}$ such that $s(v) = 0$ for every vertex $v\in V(G)$. The $null$ $set$ of a graph $G$, denoted by $N(G)$, is the set of all natural numbers $h \in \mathbb{N}$ such that $G$ admits a zero-sum h-magic labeling.

Recently, Saieed Akbari, Farhad Rahmati and Sanaz Zare $ \cite {sai} $ obtained the following interesting results about magic labeling of regular graphs.

\bigskip
\begin{itshape}

{\bf Theorem 1.1 $^{\cite {sai}}$}   \ Let $G$ be an $r$-regular graph ($r\geqslant3$, $r\neq5$). If $r$ is even, then $N(G) = \mathbb{N}$,
otherwise $\mathbb{N} \setminus \{2, 4\}$ $\subseteq N(G)$. Furthermore, if $r$ ($r\neq5$) is odd and $G$ is a 2-edge connected $r$-regular graph, then $N(G) = \mathbb{N} \setminus \{2\}$.

\end{itshape}
\bigskip

They also proposed the following conjecture in $\cite {sai}$.

\bigskip

\begin{itshape}

{\bf Conjecture}   \  Every 5-regular graph admits a zero-sum 3-magic labeling.

\end{itshape}

\bigskip

In this  paper, we give an affirmative answer to this conjecture. The following lemma is essential in the proof of the conjecture.

\bigskip

\begin{itshape}

{\bf Lemma 1.1 $^{\cite{yu}}$}   \  Let $G$ be a graph of even order with degree sequence $d$=($d_1$, $d_2$, $\dots$, $d_{n}$). If $\tilde{d}$=($d_1-1$, $d_2-1$, $\dots$, $d_{n}-1$) is also a degree sequence of some graph, then $G$ has a 1-factor.

\end{itshape}

\bigskip

More information and related references concerning magic labeling of graphs can be seen in $\cite {sai}$.

\bigskip

 \noindent {\bf 2. Main Results}

\bigskip

In this section, we will give a proof of the Conjecture.

If a graph $G$ has vertices $v_1$, $v_2$, $\dots$, $v_{n}$, the sequence $d$=($d_1$, $d_2$, $\dots$, $d_{n}$) is called the $degree$ $sequence$ of $G$, where $d_{i}=d(v_{i})$ for $i=1, 2, \dots, n$. A nonincreasing and nonnegative integer sequence $d$=($d_1$, $d_2$, $\dots$, $d_{n}$) is $graphical$ if there is a simple graph with degree sequence $d$. It is obvious that the conditions $d_{i}\leqslant n-1$  for all $i$, and $\sum _{i=1}^{n} d_{i}$ being  even are necessary for a sequence to be graphical. Firstly, the following lemma will be obtained.

\bigskip

\begin{itshape}

{\bf Lemma 2.1}   \  Let $n$ be a positive even number, and $d$=($d_1$, $d_2$, $\dots$, $d_{n}$) be a sequence of nonnegative integers. If $d_1$ = $d_2$ =  $\dots$ = $d_{n}$ = 5 and $n\geqslant 6$, or $d_1$ = $d_2$ =  $\dots$ = $d_{n}$ = 4 and $n\geqslant 6$, then $d$ is graphical.

\end{itshape}

\bigskip

{\bf Proof}  \  \ For convenience, we let $G_{n}$ denote the corresponding graph related to the sequence $d$=($d_1$, $d_2$, $\dots$, $d_{n}$).

Firstly, we prove that if $d_1$ = $d_2$ =  $\dots$ = $d_{n}$ = 5 and $n\geqslant 6$  then $d$ is graphical. The proof is by induction on $n$. If $n=6$, then it is a obvious result since the complete graph $K_6$ being the graph $G_{6}$ with degree sequence $(5, 5, 5, 5, 5, 5)$. When $n=8$, the corresponding graph $G_8$ is obtained from $G_6$ through the following construction. Let $V(G_6)$=\{$v_1$, $v_2$,  $\dots$,  $v_6$\}. Firstly, we add two new vertices $v_7$ and $v_8$ to $G_6$, and add an edge connecting $v_7$ and $v_8$. Secondly, we select, in $G_6$, two different matchings $M_1$ and $M_2$ with $M_1\cap M_2=\emptyset$ and $|M_1|$=$|M_2|$=2. Deleting the four edges in  $M_1\cup M_2$ from $G_6$, and connecting the four vertices in $M_1$ to $v_7$, the other four vertices in $M_2$ to $v_8$, we get the graph $G_8$ with degree sequence $(5, 5, 5, 5, 5, 5, 5, 5)$. Now, suppose that $n=2(k+1)\geqslant10$. By induction hypothesis the $2k$-elements $sequence$ $(5, 5, \dots, 5)$ is graphical and the corresponding graph is $G_{2k}$. So the graph $G_{2(k+1)}$ can be obtained from $G_{2k}$ through the same procedure as that of $G_6$ to $G_8$, and the proof is complete.

As for the case $d_1$ = $d_2$ =  $\dots$ = $d_{n}$ = 4 and $n\geqslant 6$, we also through the induction on $n$. If $n=6$, then it is an easy work to find a 4-regular graph $G_6$ with degree sequence $(4, 4, 4, 4, 4, 4)$. When $n=8$, the corresponding graph $G_8$ is obtained from $G_6$ through the following operation. Let $V(G_6)$=\{$v_1$, $v_2$,  $\dots$,  $v_6$\}. Firstly, we add two new vertices $v_7$ and $v_8$ to $G_6$, and select, in $G_6$, two different matchings $M_1$ and $M_2$ with $M_1\cap M_2=\emptyset$ and $|M_1|$=$|M_2|$=2. Deleting the four edges in  $M_1\cup M_2$, and connecting the four vertices in $M_1$ to $v_7$, the other four vertices in $M_2$ to $v_8$, we get the graph $G_8$ with degree sequence $(4, 4, 4, 4, 4, 4, 4, 4)$. Now, suppose that $n=2(k+1)\geqslant10$. By induction hypothesis the $2k$-elements $sequence$ $(4, 4, \dots, 4)$ is graphical and the corresponding graph is $G_{2k}$. So the graph $G_{2(k+1)}$ can be obtained from $G_{2k}$ through the same procedure as that of $G_6$ to $G_8$, and the proof is complete.     $\hspace*{\fill} \Box$

\bigskip

\begin{itshape}

{\bf Theorem 2.1}  \  Every 5-regular graph admits a zero-sum 3-magic labeling.

\end{itshape}

\bigskip

{\bf Proof}  \  \ It is obvious that every 5-regular graph $G$ is of even order since $2\cdot E(G)=5\cdot V(G)$. For $|V(G)|< 6$, the correctness of the theorem is  easily  to verify.   When $|V(G)|\geqslant 6$, according to the Lemma 1.1 and Lemma 2.1 we can get that every 5-regular graph $G$ contains a 1-factor. So, labeling the edges in the 1-factor with 2 ($\in \mathbb{Z}_{3} \setminus \{0\}$)  and the remaining edges with 1 ($\in \mathbb{Z}_{3} \setminus \{0\}$), we will get a zero-sum 3-magic labeling of the 5-regular graph. $\hspace*{\fill} \Box$

\medskip

The following theorem can be easily deduced from the Theorem 1.1 and Theorem 2.1.

\bigskip

\begin{itshape}

{\bf Theorem 2.2 }   \ Let $G$ be an $r$-regular graph with $r\geqslant3$. If $r$ is even, then $N(G) = \mathbb{N}$, otherwise $\mathbb{N} \setminus \{2, 4\}$ $\subseteq N(G)$. Furthermore, if $r$  is odd and $G$ is a 2-edge connected $r$-regular graph, then $N(G) = \mathbb{N} \setminus \{2\}$.

\end{itshape}

\bigskip

\hspace{-4mm}$\bf{Acknowledgements}$
 \ \ The authors  thank the referees for their careful reading of the
paper, and for their valuable comments.

\medskip


\end{document}